\documentclass[11pt]{amsart}

\usepackage{amsfonts} 
\usepackage{graphicx}
\usepackage{color}
\usepackage{amsthm}
\usepackage{hyperref}
\usepackage{todonotes}
\usepackage{natbib}

\usepackage{pinlabel}

\newtheorem{theorem}{Theorem}[section]

\newtheorem{corollary}[theorem]{Corollary}
\theoremstyle{definition}
\newtheorem{definition}[theorem]{Definition}

\newtheorem{example}[theorem]{Example}

\theoremstyle{remark}
\newtheorem{remark}[theorem]{Remark}

\def\a{\alpha}

\def\b{\beta}

\def\g{\gamma}

\def\G{\Gamma}

\newcommand{\nc}{\newcommand}

\nc{\p}[1]{\medskip\noindent{\em #1.}}
\nc{\margin}[1]{\marginpar{\scriptsize #1}}

\begin{document}


\title{MICC: A Tool for Computing Short Distances in the Curve Complex}

\author{Paul Glenn}
\address{Department of Mathematics \\ University at Buffalo---SUNY}
\email{paulglen@buffalo.edu}

\author{William W. Menasco}
\address{Department of Mathematics \\ University at Buffalo---SUNY}
\email{menasco@buffalo.edu}

\author{Kayla Morrell}
\address{Department of Mathematics \\ Buffalo State College---SUNY}
\email{morrelke01@mail.buffalostate.edu}

\author{Matthew Morse}
\address{Department of Mathematics \\ University at Buffalo---SUNY}
\email{mjmorse@buffalo.edu}

\begin{abstract}
The complex of curves $\mathcal{C}(S_g)$ of a closed orientable surface of genus $g \geq 2$ is the simplicial
complex whose vertices, $\mathcal{C}^0(S_g)$, are isotopy classes of essential simple closed curves in $S_g$.  Two vertices
co-bound an edge of the $1$-skeleton, $\mathcal{C}^1(S_g)$, if there are disjoint representatives in $S_g$.  A metric is obtained on $\mathcal{C}^0(S_g)$ by
assigning unit length to each edge of $\mathcal{C}^1(S_g)$.  Thus, the distance between two vertices, $d(v,w)$, corresponds
to the length of a geodesic---a shortest edge-path between $v$ and $w$ in $\mathcal{C}^1 (S_g)$.  In \cite{[BMM]},
 Birman, Margalit and the second author introduced the concept of {\em efficient geodesics} in
$\mathcal{C}^1(S_g)$ and used them to give 
a new algorithm for computing the distance between vertices.  In this note, we introduce the software package MICC ({\em Metric in the Curve Complex}), a partial implementation
of the efficient geodesic algorithm.  We discuss the mathematics underlying MICC and give applications.  In particular,
up to an action of an element of the mapping class group,
we give a calculation which produces all distance $4$ vertex pairs for $g=2$ that intersect $12$ times, the minimal number of
intersections needed for this distance and genus.
\end{abstract}

\maketitle

\section{Introduction}
\label{sec:intro}

Let $S$ or $S_g$ denote a compact, connected, orientable surface of genus $g$, where $g \geq 2$. 
A simple closed curve on $S$ is {\em essential} if does not bound a disk in $S$.
The complex of curves, introduced by Harvey \cite{[Ha]}, is the simplicial complex, $\mathcal{C}(S)$,
whose vertices (or $0$-skeleton), $\mathcal{C}^0(S)$, are isotopy classes of essential simple closed curves; and,
whose edges of the $1$-skeleton, $\mathcal{C}^1(S)$, connect vertices that have disjoint representatives.
For the remainder of this note, ``curve'' will mean ``simple closed curve''.
By declaring that each edge of $\mathcal{C}^1(S)$ has length $1$, we endow $\mathcal{C}^0(S)$
with a metric.  Specifically, an {\em edge path} is a sequence of vertices $\{v=v_0 , v_1 , \cdots , v_n=w\} $ such that $d(v_i , v_{i+1})=1$. 
A {\em geodesic path} joining $v$ and $w$ is a shortest edge-path.  The {\em distance}, $d(v,w)$, between arbitrary vertices is the length of a geodesic path.  Since it is known that the complex of curves is connected, which was stated by Harvey \cite{[Ha]} and followed from a previous argument of Lickorish \cite{[Li]}, the value $d(v,w)$ is well-defined for all vertex pairs.

We note that if $d(v,w)=2$, there is a vertex $\bar\g \in \mathcal{C}^0(C)$ and curve representatives in $S$,
$\alpha \in v$, $\beta \in w$ and $\gamma \in \bar\g$, such that $\alpha \cap \beta \not= \emptyset$
and $\gamma \subset S \setminus (\alpha \cup \beta)$.  
The generic situation  (when some component of
$S \setminus (\alpha \cap \beta)$ has Euler characteristic less than zero) is that there are infinitely many isotopically distinct choices for $\gamma \subset S \setminus (\alpha \cap \beta)$ and, thus, infinitely many possible geodesics for distance $2$. 
In this case, the existence of infinitely many geodesics at distance $2$ forces infinitely many geodesics for all distances.  
It is this infinite local pathology which makes finding an effective distance computing algorithm challenging.

The curve complex was first introduced by Harvey \cite{[Ha]}.  Its coarse geometric properties were first studied extensively by Masur--Minsky \cite{[MM1], [MM2]}.
The complex of curves has proved a
useful tool for the study of hyperbolic $3$-manifolds, mapping class groups and Teichm\"uller theory.  In particular, in \cite{[MM2]}, it was established that Teichm\"uller space is
quasi-isometric to the complex of curves and is therefore $\delta$-hyperbolic.  Here, $\delta$-hyperbolic means that geodesic triangles in $\mathcal{C}^1(S)$ are $\delta$-thin:
any edge is contained in the $\delta$-neighborhood of the union of the other two edges.  In \cite{[A]}, Aougab established uniform hyperbolicity---$\delta$ can be chosen independent of genus (for $g \geq 2$).  In spite of this considerable advancement in understanding the coarse geometry of the complex of curves, the development of tools intended to explicitly compute distance has been difficult.

In 2002, Jason Leasure proved  the existence of an algorithm to compute the distance between two vertices of $\mathcal{C}^0(S)$ (\cite{[L]}, Corollary 3.2.6).   Later, other algorithms were discovered by Shackleton \cite{[Sh]} and Webb \cite{[W]}, but none of these algorithms were studied seriously from the viewpoint of doing explicit computations, and all seem unsuitable for that purpose.

Recently Birman, Margalit and the second author \cite{[BMM]} have given a new algorithm---{\em the efficient geodesic algorithm}---and we
have developed an implementation of it called the {\em Metric in the Curve Complex} (MICC).  Applications of MICC we will present in this
note include:
\begin{itemize}
\item[(i)] establishing that the minimal geometric intersection number for vertices of $\mathcal{C} (S_2)$ with distance four is $12$,
\item[(ii)] listing of all vertex pairs (up to an action of an element of the mapping class group) of $\mathcal{C} (S_2)$ with distance four and having minimally positioned representatives with intersection number at most $25$, and
\item[(iii)] producing an explicit example of two vertices of $\mathcal{C} (S_3)$ that have distance
four and intersection number $29$
\end{itemize}

The key idea in \cite{[BMM]}  is the introduction of a new class of geodesics, {\em efficient geodesics}.   They are not the same as the `tight geodesics' that have dominated almost all published work on the curve complex following their introduction in \cite{[MM1],[MM2]}, however they share with tight geodesics the nice property that there are finitely many efficient geodesics between any two fixed vertices in $\mathcal{C}(S)$.

For convenience, for a pair of curves, $(\a , \b)$, we will refer to a component of
$(\a \cup \b) \setminus (\a \cap \b)$ as a {\em segment}.
We will use a slightly weaked definition for efficient geodesic than that given in \cite{[BMM]}.

\begin{definition}
\label{D: IE}
Let $v, w \in \mathcal{C}^0(S)$ with $d(v,w) \geq 3$.   An oriented path $v=v_0, \dots, v_n=w, \ n\geq 3$, in $\mathcal{C}^0(S)$  
is  {\em initially efficient} if there are representatives $\alpha_0 \in v_0$, $\alpha_1 \in v_1$ and $\alpha_n \in v_n$ such that
$| \alpha_1 \cap b| \leq n-2 $ for any segment $b \subset \alpha_n \setminus \alpha_0$. We say $v = v_0 , \cdots , v_n = w$ is {\em efficient} if $v_k , \cdots , v_n$ is initially
efficient for each $ 0 \leq k \leq n-3$ and the oriented path $v_{n} , v_{n-1} , v_{n-2} , v_{n -3}$ is also initially efficient.
\end{definition}

The efficient path algorithm is a consequence of the following.

\begin{theorem} {\rm (Theorem 1.1 of \cite{[BMM]})}
\label{T:simplify}
Let $ g \geq 2$, and let $v$ and $w$ be two vertices of $\mathcal{C} (S_g)$ with $d(v,w) \geq 3$.
There exists an efficient geodesic from $v$ to $w$, and in fact there are finitely many.
\end{theorem}


When $n=3$, notice that an efficient geodesic $v=v_0 , v_1 , v_2 , v_3=w$
yields an oppositely oriented efficient geodesic, $w=v_3 , v_2 , v_1 , v_0=v$.  That is, distance $3$ vertices have non-oriented
efficient geodesics. 
Thus, for corresponding representatives $\a_0, \a_1, \a_2, \a_3$, we have that $\a_1$ (respectively $\a_2$) will intersect
segments of $\a_3 \setminus \a_0$ (respectively of $\a_0 \setminus \a_3$) at most once.
From this observation, we will establish the following test for distance $\geq 4$ which MICC implements.

\begin{theorem}
\label{theorem: BMM test}{\rm({\bf Distance $\geq 4$ Test})}
Let $v, w$ be vertices with $d(v,w) \geq 3$.  Let $\G\subset \mathcal{C}^0(S)$
be the collection of all vertices such that the following hold:
\begin{enumerate}
\item for $\bar\g \in \G$, we have $d(v,\bar\g)=1$; and
\item for $\bar\g \in \G$, there exist representatives $\alpha, \beta,\gamma$ of $v,w,\bar\g$ respectively, such that for
each segment $b\subset \beta\setminus \alpha$ we have $| \gamma\cap b | \leq 1$.
\end{enumerate} 
Then $d(v,w) \geq 4$ if and only if $d(\bar\g,w) \geq 3$ for all $\bar\g \in \G$.  Moreover, the collection $\G$ is finite.
\end{theorem}

 \begin{remark}
 \label{remark on Z}
Keeping with our previous observation regarding non-oriented efficient geodesics at distance $3$, we can
flip the roles of $v$ and $w$.  Thus, the test can also
be stated in terms of $d(v,\bar\g^\prime) \geq 3$ for $\bar\g^\prime \in \G^\prime$ where:
\begin{itemize}
\item[1.] If $\bar\g^\prime \in \G^\prime$ then $d(\bar\g^\prime,w)=1$.
\item[2.] If $\bar\g^\prime \in \G^\prime$, there exists representatives, $\alpha \in v$, $\beta \in w$ and $\g^\prime \in \bar\g^\prime$ such that
for all segments, $a \subset \alpha \setminus \beta$, we have $| a \cap \gamma^\prime | \leq 1$.
\end{itemize}
 \end{remark}
 
 Before giving the proof of the Distance $\geq 4$ Test, we recall a useful concept
 and its implications.
 Let $\alpha ,\beta \subset S$ be a pair of curves such that $|\alpha \cap \beta|$ is minimal with respect to isotopies of $\beta$.
 That is, $\alpha$ and $\beta$ are {\em minimally positioned}.  Determining when $\a$ and $\b$ are minimally positioned
 is straightforward due to the bigon criterion (Propostion 1.7, \cite{[FM]})---no disc component of $S \setminus (\a \cup \b)$
 has exactly two segments of $(\a \cup \b) \setminus (\a \cap \b)$ in its boundary.
 We say $\alpha$ and $\beta$ (or $(\alpha,\beta)$) is a {\em filling pair} if
 $S \setminus (\alpha \cup \beta)$ is a collection of $2$-discs.  It is readily seen
 that a pair is filling on $S$ if and only if their corresponding vertices
 in $\mathcal{C}^0(S)$ are at least distance $3$ apart.
 When a minimally positioned pair of curves is not filling but still intersects, some component of
 $S \setminus (\alpha \cup \beta)$ contains an essential curve.  Thus, the corresponding vertices are distance $2$ apart.
 Algorithmically determining whether a minimally positioned pair is filling, or not, requires simple tools coming from classical topology.
For $\alpha \cap \beta \not= \emptyset$, let $N(\alpha \cup \beta) \subset S_g$ be a regular neighborhood.
The genus of this neighborhood, $genus(N(\alpha \cup \beta))$, can be algorithmically computed as discussed in \S\ref{sec: genus}.
(For an oriented surface $\Sigma$ with boundary, recall $genus(\Sigma) = 1 - \frac{ \chi(\Sigma) + |\partial \Sigma|} {2}$,
where $\chi(\Sigma)$ is the Euler characteristic and $\chi (N(\alpha \cup \beta)) = -|\a \cap \b|$.)
If $genus(N(\alpha \cup \beta)) < g$, then
a component of $S_g \setminus (\alpha \cup \beta)$ contains an essential curve of $S_g$ and the vertices that
$\alpha$ and $\beta$ represent are distance $2$ apart.  If $genus(N(\alpha \cup \beta))=g$, those vertices
are distance at least $3$ apart.  We will see in \S~\ref{sec: genus} that this {\em filling calculation} can be readily implemented.  However,
if one is handed a nice enough presentation of $\a$ and $\b$ in $S$, determining whether they are a filling pair can be done by inspection.
For example, we will do such filling determinations in Example~\ref{example:Hempel}.

We now give the proof of Distance $\geq 4$ Test.

\noindent {\bf Proof:}  From the above discussion we see that the assumption, $d(v,w) \geq 3$, translates into
considering only minimally positioned filling pairs in $S$.  To determine whether the associated vertices
in $\mathcal{C}^0(S)$ of a filling pair are at distance $\geq 4$, we need only determine that they are not
at distance 3. 

Thus, suppose $\alpha$ and $\beta$ represent classes $v$ and $w$ such that $d(v,w) \geq 3$.  Assume there exists a length $3$ path
$v =v_0, v_1, v_2, v_3=w$.  From Theorem \ref{T:simplify}, we can further assume that this path is initially efficient.
In particular, for representative $\alpha = \a_0 , \a_1 , \a_2, \a_3=\beta$ of the vertices of this path, respectively,
we can assume $\alpha_1$ intersects segments of $\alpha_3 \setminus \alpha_0$ at most once.
Thus, $v_1$ is an element of the set $\G$.  But, since $d(v_1 , v_3(=w)) =2$, we need only establish that
$d(\bar\g,w)\geq 3$ for all $\bar\g \in \G$ to contradict the assumption that there was a distance $3$ path.

The fact that the set $\G$ is finite is due to $S \setminus (\alpha_0 \cup \alpha_3)$
being a collection of $2$-discs and representatives of any vertex of $\G$ having bounded intersection with any such
$2$-disc component.  The stated test for $\geq 4$ follows.
\qed

\begin{example}  \label{example:Hempel}
We consider an example of a pair of curves, $\alpha$ and $\beta$, on a genus $2$ surface which represent classes that
are distance $4$ apart.  (See Figure~\ref{fig: Hempel example}.)  This is an example of J. Hempel and appears in the lecture notes \cite{[Sc]}.  These notes assert
distance $4$ for the pair without proof.  As an application of the Distance $ \geq 4$ Test, we now give a proof
establishing distance $4$ for Hempel's example.

\begin{figure}[htbp!]
\labellist
\small\hair 2pt
\pinlabel $\mathrm{arcs \  of}\ \alpha$ at 140 100
\pinlabel $\beta$ at 165 119
\pinlabel $\beta$ at 363 244
\pinlabel $\gamma$ at 97 280
\pinlabel $\gamma^\prime$ at 380 276
\endlabellist
\centering{\includegraphics[width=1.1\textwidth]{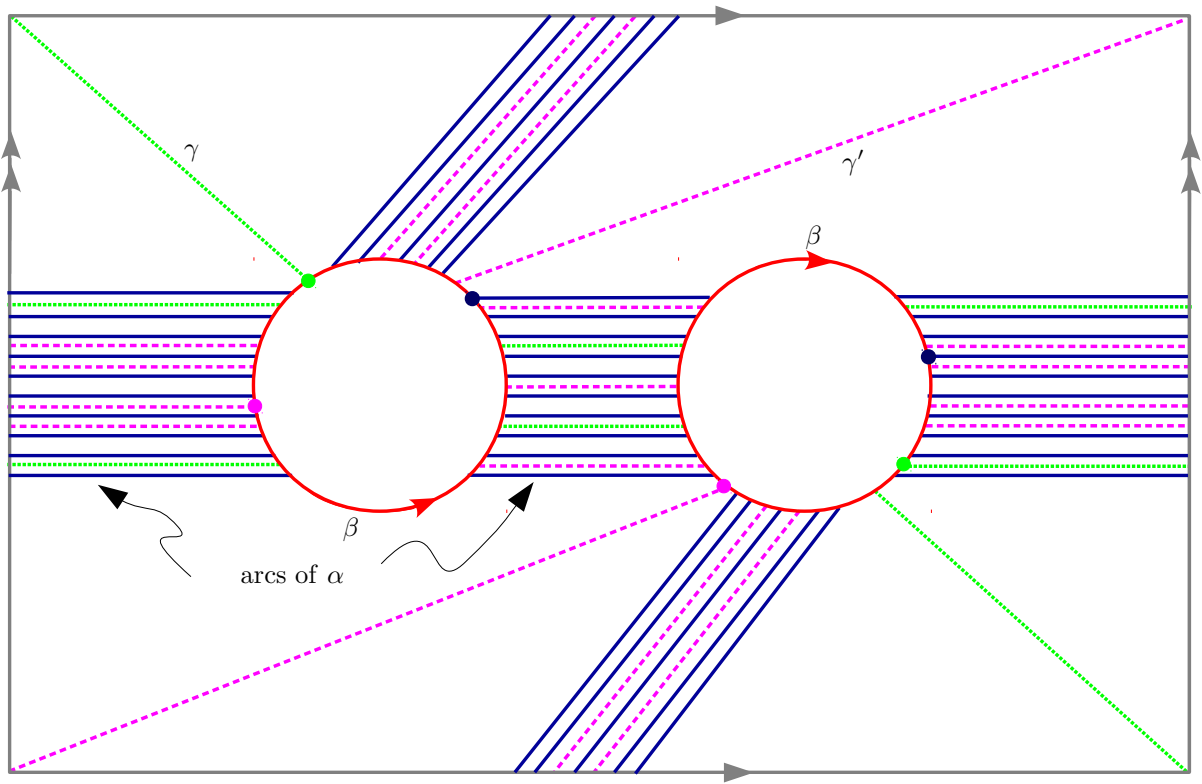}}

\caption{{\small An example due to J. Hempel} }
\label{fig: Hempel example}
\end{figure}

In Figure \ref{fig: Hempel example}, the surface $S_2$ is represented as a rectangular region minus two discs.
The gray  sides of the rectangle are identified, left-to-right and top-to-bottom, to form a torus
minus two discs.  The genus $2$ surface is obtained by identifying the two oriented red boundary curves, and
the resulting single curve is $\beta$  (The identification is initiated by lining up the six colored dots on the $\a, \g ,\g'$ curves).
These identifications induce identifications of the endpoints of the dark blue arcs, so as to form the curve $\alpha$.
By inspection, one can see that $(\a, \b)$ is a minimally positioned filling pair.

We now apply the Distance $\geq 4$ Test.  We wish to find curves, $\g$, that represent vertices, $\bar\g \in \G$.
Such a $\g$ will be in the complement of $\alpha$, intersecting any segment of $\beta \setminus \alpha$
at most once.  Three such $\g{\rm 's}$ can be immediately identified.
This is because the complement of $\alpha \cap \beta$ in $S_2$ is a collection of some 
number of $4$-gon regions and one single $12$-gon region.
The boundary of any one of these regions is an alternating joining of segments in $\alpha$ and $\beta$.  Thus,
any $4$-gon boundary has two segments in $\beta$; and, the boundary of the single $12$-gon has
$6$ segments in $\beta$.  Requiring that any $\g$ intersect segments of $\beta$ at most once forces it to either not intersect a $4$-gon, or
intersect each of the two $\beta$ segments of a $4$-gon once.  However, there six different ways a $\g$ can exit/enter
the $12$-gon, giving us three possible $\g{\rm 's}$ that intersect the $12$-gon region once.
In Figure \ref{fig: Hempel example}, the dashed green and purple curves $\g$ and $\g^\prime$ illustrate
two of the three curves generated by the exit/enter possibilities.  It is readily apparent that both
$S \setminus (\beta \cup \g)$ and $S \setminus (\beta \cup \g^\prime)$ is a collection
of $2$-discs, none of which are bigons.  Thus, the corresponding vertex pairs are at least distance $3$ apart.
The remaining possibilities for a $\gamma$ can be dealt with
in a similar straightforward manner (theoretically, there are also $\gamma{\rm 's}$ that intersect the $12$-gon region
$2$ and $3$ times). $\diamond$
\end{example}

The $\G$-calculation above illustrates the primary computing
capabilities of the MICC
software package \cite{[MICC]}.   MICC is a computational
tool that can determine whether the distance between two vertices in $\mathcal{C}(S_{g\geq2})$ is $2$, $3$,  or $\geq 4$.  Its input is readily produced from any representation of two curves on a closed surface.
It has functionality that can be used to search for new curve pairs or manipulate existing examples.  Its output
can be used to construct geodesic paths between curves of short distances.

As such, MICC is an additional tool scholars can utilize in answering a number of basic questions about
the local pathology of the complex of curves.  As an illustration, we consider the relationship between distance and minimal intersection number.
It is known that the theoretical minimal intersection
number for a filling pair on a $S_g$ is $2g-1$ due to the Euler characteristic of the surface.  For $g=2$, this theoretical minimum is not realizable and the
realizable minimum is in fact $4$.  Recent work of Aougab and Huang \cite{[AH]} has given a construction
for realizing the theoretical minimum for $g\geq3$.  Additionally, they show that all such minimum filling pairs are
distance $3$.  For fixed $g \geq 2$, using his uniform hyperbolicity result, Aougab
proved that the theoretical minimum intersection number grows exponentially as a function of distance (Theorem 1.2, \cite{[A]}).
Also, Aougab and Taylor \cite{[AT]} give a recipe for producing filling pairs at a given distance whose intersection numbers are close to the minimum in an asymptotic sense; see their paper for the precise statement.  Ido, Jang and Kobayashi \cite{[IJK]} also
have a construction for producing filling pairs of a prescribed distance.
The arguments in these last three citations employs the high power machinery of Masur and Minsky, including the {\em Bounded geodesic image theorem} \cite{[MM2]}.  Thus, the growth bounds and constructed examples inherit a ``coarse geometry'' quality,
which so far in the literature has not been used to produce the
exact minimal intersection number with accompanying filling pairs for a specified distance and genus.


In contrast, MICC can be used to find explicitly all minimum intersecting filling pairs of distance $4$.
Using the Distance $\geq 4$ Test, we give a ``proof of concept'' calculation that constructs all minimum
intersecting distance $4$ filling pairs in $\mathcal{C}(S_2)$ up to homeomorphism.  We can next use MICC to calculate distance for curve pairs of increasing intersection number starting at this minimum.  The result of this calculation is the following theorem.

\begin{theorem}
\label{theorem: DWH-G2-D4-I12}
The minimal intersection for a filling pair, $\alpha , \beta \subset S_2$, representing
vertices $v, w \subset \mathcal{C}^0(S_2)$, respectively, with $d(v,w)=4$ is $12$.
\end{theorem}

Combining this theorem with the Distance $\geq 4$ Test we obtain a partial test for distance four.

\begin{corollary}
\label{corollary: DWH-G2-D4-I12}
Let $d(v,w) \geq 4$ for two vertices in $\mathcal{C}^0(S_2)$.
Let $\alpha, \beta \subset S_2$ be curves in minimal position representing $v$ and $w$, respectively.
Let $\gamma \subset S_2 $ be a curve.
If 
\begin{enumerate}
\item $\gamma \cap \beta = \emptyset$ and $|\gamma \cap \alpha| < 12$, or
\item $\gamma \cap \alpha = \emptyset$ and $|\gamma \cap\beta|<12$,
\end{enumerate}
 then $d(v,w)=4$.
\end{corollary}

The proof of concept calculation for Theorem \ref{theorem: DWH-G2-D4-I12} involves finding all solutions to
an integer linear programing problem so as to identify all potential candidates for minimally intersecting
distance four curve pairs.  However, such a comprehensive search is not necessarily needed to find examples
of distance $\geq 4$ pairs.
Utilizing all of the functionality of MICC, one can ``experiment'' with different curve pairs in a search distance four pairs.
Remark \ref{remark: experiment} discusses how such experimentation led to the discovery of the first known
filling pair representing distance $4$ vertices in $\mathcal{C}^0(S_3)$.
In particular, we have the following result:

\begin{theorem}
\label{theorem: DWH-G3-D4-I29}
The minimal intersection number for a pair of filling curves in $S_3$ that represent distance $4$ vertices in $\mathcal{C}^0(S_3)$
is less than or equal to $29$.
\end{theorem}


The outline for our paper is as follows.
In \S\ref{sec: DWH}, we discuss a method of representing
a filling pair on $S_{g \geq 2}$.
In \S\ref{sec: DWH-G2-D4-I12}, we give the proof of Theorem \ref{theorem: DWH-G2-D4-I12}.  In particular, the proof
can be viewed as giving a general strategy for calculating theoretical minimal intersections of distance $4$ filling pairs
for any higher genus.  This strategy employed in a limited manner allowed us to verify
that the explicit example given in \cite{[BMM]} (cf. \S2) of a distance $4$ vertex pair in $\mathcal{C}^0(S_2)$
establishes Theorem \ref{theorem: DWH-G2-D4-I12}.  We finish \S\ref{sec: DWH-G2-D4-I12}
with an analysis of a genus $3$ distance $4$ pair establishing Theorem \ref{theorem: DWH-G3-D4-I29}.
Finally, in \S\ref{sec: MICC commands}, we discuss the
complete functionality of MICC.  To make this discussion concrete, we illustrate the range of MICC commands with a running
example.

At the end of this manuscript, we attach the current known {\em spectrum} of pairs of distance $4$ or greater in $\mathcal{C}^0(S_2)$ with up to $14$ intersections.  The full distance $4$ or greater spectrum in $\mathcal{C}^0(S_2)$ with
up to $25$ intersections is available at \cite{[MICC]}.

\noindent
{\em An expository remark}---Throughout we will continue to use $\a$ and $\b$ as representatives of $v$ and $w$,
respectively.  Similarly, indexed $\a_i$ curves will be used as representatives of indexed $v_i$ vertices.
This is meant to be consistent with the notation used in \cite{[BMM]}.  For all other curves in the surface, we will use
various ``flavors'' of $\gamma$, and $\bar\g$ will denote the corresponding vertex.  We will always assume that
any pairing of curves are in minimal position.

\section{Representations of pairs of curves.}
\label{sec: DWH}
\subsection{Disc with handles}
\label{subsec: DWH}
In deciding how to represent curves on surfaces, we must first choose how we will represent closed oriented surfaces.
Representing surfaces with boundary as {\em disc with handles} (or $DWH$, example in Figure~\ref{fig: DWH-G2-2}) is a well known method among working
geometers and topologists, and it can be readily adapted to closed surfaces in our situation.  For an essential
curve $\a \in S_{g \geq 2}$, we can split $S$ along $\alpha$ to produce a surface $\hat{S}$ having two boundary
curves, $\partial_+$ and $\partial_-$.  If $\alpha$ is separating, $\hat{S}$ will have two connected components, i.e. $\hat{S} = \hat{S}^1 \cup \hat{S}^2$.  The genus of each component
will be less than $g$ and their sum would be ${\rm genus}(\hat{S}^1) + {\rm genus}(\hat{S}^2) = g$.  As such, a $DWH$ representation
will have $\hat{S}^1$ (respectively, $\hat{S}^2$) being a single $2$-disc with $2\times{\rm genus}(\hat{S}^1)$
(respectively, $2\times{\rm genus}(\hat{S}^2)$) $1$-handles attached. We recover $S$ by gluing the boundary of these two components together.

When $\alpha$ is non-separating $\hat{S}$ will be a $g-1$ connected surface with two boundary curves $\partial_+$ and $\partial_-$.
A $DWH$ representation of $\hat{S}$ would
be a single $2$-disc with $2 \times {\rm genus}(\hat{S}) + 1$ $1$-handles attached.  As before, $S$ is recovered
by giving a gluing of its two boundary curves.

Now suppose $(\alpha,\beta)$ is a filling pair in minimal position on $S$.
We obtain a $DWH$ representation of $S$ by splitting along $\alpha$.
Since $\alpha$ and $\beta$ are in minimal position, we know that $\beta \cap \hat{S}$ will be a collection of properly embedded essential arcs,
$\{\omega_1 , \cdots , \omega_k \}$.  We require that each $1$-handle in our $DWH$ representation contain at least one $\omega$-arc. 
For example, skipping ahead to Figure \ref{fig: DWH-G2-2}, the arcs with labels $w_1, \ w_2, \ w_3$ are the needed arcs.
We refer to such an arc as a {\em co-cores} of a $1$-handle.
Skipping ahead to Figure \ref{fig: DWH-G2-D4-I12}, the reader will find an example of
a genus $1$ $DWH$ surface with two boundary components.
The properly embedded essential black arcs are examples of $\omega$-arcs and each $1$-handle contains at least one such arc. 

Since we must be able to recover both $\beta$ and $S$ by a gluing of $\partial_+$ and $\partial_-$, we must have
$|(\cup_{1\leq j \leq k } \ \omega_j ) \cap \partial_+ | = |(\cup_{1\leq j \leq k }\  \omega_j ) \cap \partial_- |$.
More precisely, for $\partial \omega = p_1 \cup p_2$
there are three possible configurations:  $\omega$ is a $++$ arc (respectively,
$--$ arc) when $p_1 \cup p_2 \subset \partial_+$ (respectively, $p_1 \cup p_2 \in \partial_-$); and, $\omega$ is a $+-$ arc when
$p_1 \in \partial_+$ and $p_2 \in \partial_-$.  Thus we can have any
number of $+-$ $\omega$-arcs but the number of $++$ arcs must equal the number of $--$ arcs.

Finally, observe that the $\omega$-arcs divide both $\partial_+$ and $\partial_-$ into $k$ intervals.  In order to specify a gluing of
$\partial_+$ and $\partial_-$, we orient and cyclically label these $k$ intervals, by convention, $0$ through $k-1$.
Again, referring to
Figure \ref{fig: DWH-G2-D4-I12}, we illustrate a gluing of $\partial_+$ and $\partial_-$ by the $0$ through $11$ labels, shown in red and blue respectively.

\subsection{Strategy for constructing examples}
\label{subset: constructing examples}
\begin{figure}[htbp]
\labellist
\small\hair 2pt
\pinlabel $w_6$ at 116 100
\pinlabel $w_2$ at 107 245
\pinlabel $w_4$ at 340 123
\pinlabel $w_1$ at 73 305
\pinlabel $w_3$ at 380 305
\pinlabel $w_5$ at 225 90
\pinlabel $w_6$ at 625 93
\pinlabel $w_2$ at 619 245
\pinlabel $w_4$ at 840 123
\pinlabel $w_1$ at 585 305
\pinlabel $w_3$ at 885 305
\pinlabel $w_5$ at 735 94
\endlabellist
	\centering
		\includegraphics[width=0.90\textwidth]{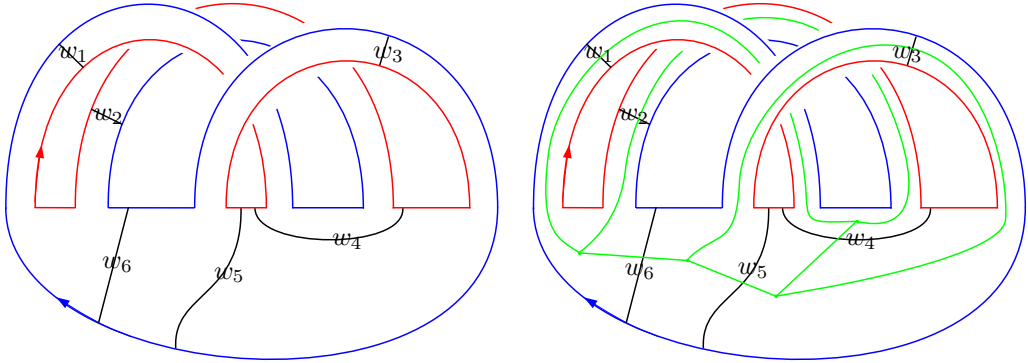}
	\caption{{\small  The left illustration is a genus one surface with two boundary curves---coded red and blue.  $C$ is a maximal
	collection of $6$ weighted arcs.  The weights, $w_1 , w_2 , w_3 , w_4 , w_5 , w_6$, are non-negative integers.  The green graph
	in the right illustration is $G(C)$, the dual graph.  Each edge of $G(C)$ intersects exactly one arc of $C$ once.}}
	\label{fig: DWH-G2-2}
\end{figure}

We can reverse engineer this construction with an eye towards finding filling pairs of distance greater than $3$.  Suppose we are interested
in finding such a filling pair $(\alpha,\beta)$ with $\alpha$ non-separating.  Any associated $\hat{S}$ will be a connected $DWH$ with two boundary curves.
Initially, let us fill $\hat{S}$ with a maximal collection, $A$, of properly embedded essential arcs that are pairwise non-parallel.  We specify
$2g +1$ arcs to be arcs associated to the $1$-handles.  Thus, when we split the $DWH$ along these $2g+1$ arcs, we obtain the underlying
$2$-disc.  Figure \ref{fig: DWH-G2-2}-left illustrates such a configuration for $g=2$.

To obtain the collection of $\omega$-arcs that will be used to produce a curve $\beta$, we will assign weights to each arc of $A$ such that there
is a reasonable expectation that the Distance $\geq 4$ Test is satisfied.  Once we have determined the weight of an arc of $A$, we will place that number of parallel copies
of the arc in the $DWH$.  To this end, we consider the dual graph, $G(A)$, to $A$ in $\hat{S}$.  (Figure
\ref{fig: DWH-G2-2}-right illustrates such a dual graph for the left configuration.)  In graph-theoretic terms, we consider \textit{elementary
circuits}---edge paths that form simple loops---in $G(A)$.
Each elementary circuit, $\gamma \subset G(A)$, represents a possible vertex in $\G$ of the Distance $\geq 4$ Test.  Let
$\Sigma(\gamma)$ be the sum of the weights of all the arcs in $A$ that a circuit $\gamma$ intersects.
Since our speculative $\beta$ curve will be the union of
arcs parallel to those in $A$ and $(\gamma,\beta)$ should be a filling pair for any circuit $\gamma$, we require that $\Sigma(\g)$ be greater than or equal to the
minimal intersection number for a filling pair, i.e. $4$ when $g=2$ and $2g -1$ when $g>2$ \cite{[AH]}.  Thus, if $\{ \g_1 , \cdots , \g_m \} \subset G(A)$ is the complete set
of elementary circuits, for each circuit we get an inequality of the form $\Sigma(\g_i) \geq 2g -1 $ (when $g>2$) or $\geq 4$ (for $g=2$).  This gives us $m$ inequalities that make
up an {\em integer linear program} (ILP).
We add to this the equality that states the sum of the weights of $++$ arcs equals the sum of the weights of
$--$ arcs.  These combined equations are the constraints for the object equation, $P$, the sum of all the weights which we wish to minimize.

Figure \ref{fig: DWH-G2-2}-right illustrates this correspondence between a $G(A)$ and a ILP when $g = 2$.  In particular,
for this dual graph, there are $6$ elementary circuits which yield $6$ weight equations.  We get a seventh equation coming
from having the weights of the $++$ and $--$ arcs being equal.
All this yields the following constraints for minimizing $P = w_1 + w_2 + w_3 + w_4 + w_5 + w_6 $:
\begin{eqnarray}
\label{LLP}
\begin{tabular}{r c r c r c r c r c r}
$w_1$ & $+$ & $w_4$ & $+$ & $w_5$ & $+$ & $w_6$ & $\geq$ & $4$ \\
$ w_2$ & $+$ & $w_4$ & $+$ & $w_5$ & $+$ &$w_6$ & $\geq$ & $4$ \\
$\ $       & $ \ $& $\ $    & $ \ $ & $w_3$ & $+$ & $w_5$ & $\geq$ & $4$ \\
$ \ $ & $\ $ & $ \ $ & $\  $ & $w_1$ &$+$ & $w_2$ & $\geq$ & $4$ \\
$ w_1$ & $+$ & $w_3$ & $+$ & $w_4$ & $+$ & $w_6$ & $\geq$ & $4$ \\
$ w_2$ & $+$ & $w_3$ & $+$ & $w_4$ &$+$ & $w_6$ & $\geq$ & $4$ \\
$\ $       & $ \ $& $\ $    & $ \ $ & $ \ $ & $ \ $ & $w_4$ & $=$ & $w_6$ \\
$\ $       & $w_1 ,$& $w_2 , $    & $w_3 , $ & $w_4 , $ & $w_5 ,$ & $w_6$ & $\geq$ & $0$
\end{tabular}
\end{eqnarray}
Using any popular computing software (e.g. Maple) one can readily check that the optimal value of this $P$ is $8$.
Thus, we have a theoretical minimal intersection number for a filling pair of distance $4$ on $S_2$.
This optimal value happens to be uniquely realized by the {\em weight solution} $[w_1 , w_2 , w_3 , w_4 , w_5 , w_6] = [2,2,2,0,2,0]$.
Placing the corresponding set of $8$ $\omega$-arcs into the $DWH$ of Figure  \ref{fig: DWH-G2-2}-left, we then have $8$ possible
ways to identify $\partial_+$ to $\partial_-$.  As discussed in \S\ref{subsec: perm}, MICC's permutation functionality can be
used to check which of these boundary identifications will result in a single $w$ curve.  For $\omega$-arcs corresponding to
the weight solution $[2,2,2,0,2,0]$ there happen to be four identifications that yield a $(\alpha,\beta)$ filling pair.  Finally, by employing
MICC's distance functionality (described in  \S\ref{subsec: distance}) to search for all elementary circuits, one can determine
that all of these intersection $8$ filling pairs whose corresponding vertices have distance $3$ in $\mathcal{C}^0(S_2)$.


\section{Proofs of main results.}
\label{sec: DWH-G2-D4-I12}
Although the case counting is extensive, the discussion in \S\ref{sec: DWH} gives us a straight forward strategy for
proving Theorem \ref{theorem: DWH-G2-D4-I12}.
\subsection{Proof of Theorem \ref{theorem: DWH-G2-D4-I12}.}
\begin{figure}[htbp]
\labellist
\small\hair 2pt
\pinlabel $w_3$ at 163 100
\pinlabel $w_2$ at 245 345
\pinlabel $w_1$ at 120 330
\pinlabel $w_3$ at 666 100
\pinlabel $w_2$ at 748 347
\pinlabel $w_1$ at 623 330
\endlabellist

	\centering
		\includegraphics[width=0.90\textwidth]{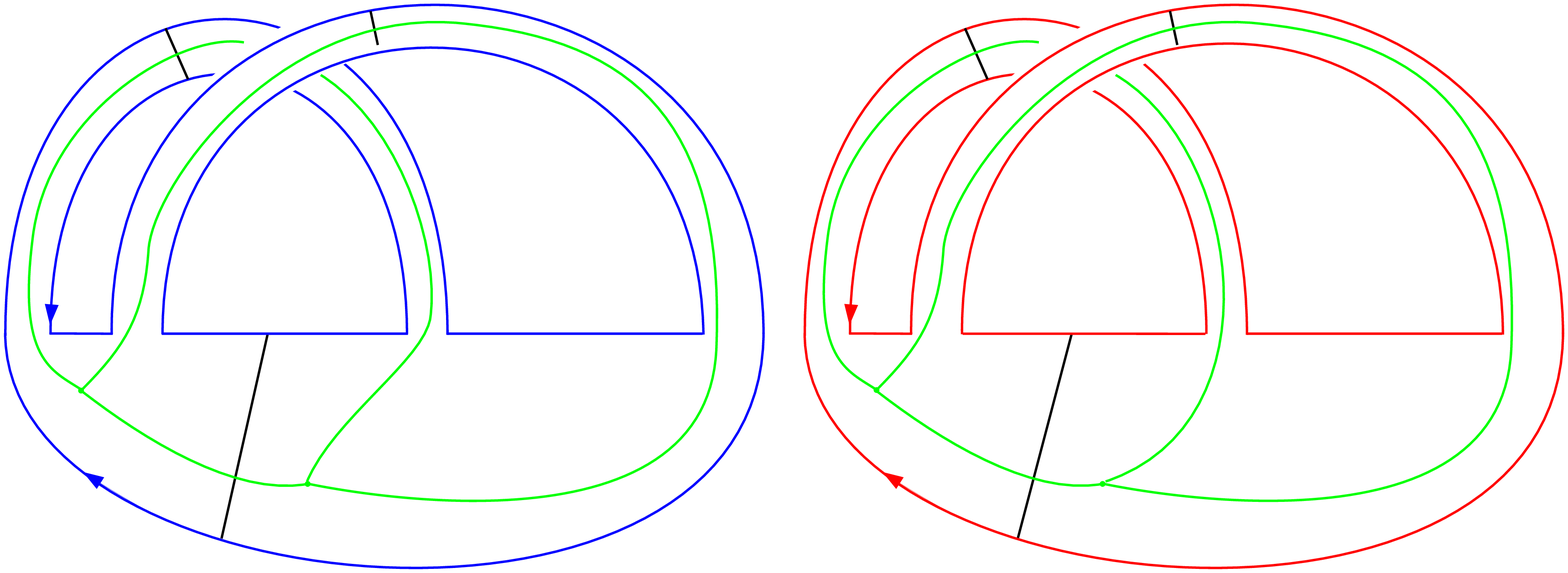}
	\caption{{\small  }}
	\label{fig: DWH-G2-4}
\end{figure}

First, we need to generate all possible genus $2$ $DWH$ diagrams with weighted $\omega$-arcs so that we may then generate the corresponding ILP's.  Initially we divide this generating process into two cases corresponding
to whether $\alpha$ is a separating or non-separating curve.  
If $\alpha$ is separating, we have exactly one possible $DWH$ diagram with weighted $\omega$-arcs.  $\alpha$ splits $S_2$ into two
genus one surfaces with boundary, and each genus one surface has three weighted $\omega$-arcs.  (See Figure \ref{fig: DWH-G2-4}.)
We leave it to the reader to generate the unique ILP in this situation.

For the ILP corresponding to Figure \ref{fig: DWH-G2-4},
we determined that the theoretical minimal intersection number for a filling pair $(\a,\b)$
of distance $4$ on $S_2$ is $12$; and, intersection $12$ is uniquely realized when all the weights equal $2$.
IBM's software package {\em CPLEX Optimization Studio} \cite{[I]} was utilized in solving this and all other ILP's
in this paper.

Next, employing MICC's permutation functionality, we determined that there are six possible identifications of $\partial_+$ and $\partial_-$
which yield a single $\b$ curve.   However, MICC's distance functionality determined that all of these filling pairs of intersection $12$ were
distance $3$.  So we turn to the non-separating case for $\alpha$.
\begin{figure}[htbp]
\labellist
\small\hair 2pt
\pinlabel $w_6$ at 106 110
\pinlabel $w_2$ at 245 349
\pinlabel $w_4$ at 315 135
\pinlabel $w_1$ at 120 332
\pinlabel $w_3$ at 297 280
\pinlabel $w_5$ at 255 85
\pinlabel $w_6$ at 606 110
\pinlabel $w_2$ at 747 349
\pinlabel $w_4$ at 855 124
\pinlabel $w_1$ at 620 332
\pinlabel $w_3$ at 797 280
\pinlabel $w_5$ at 755 85

\endlabellist

	\centering
		\includegraphics[width=0.90\textwidth]{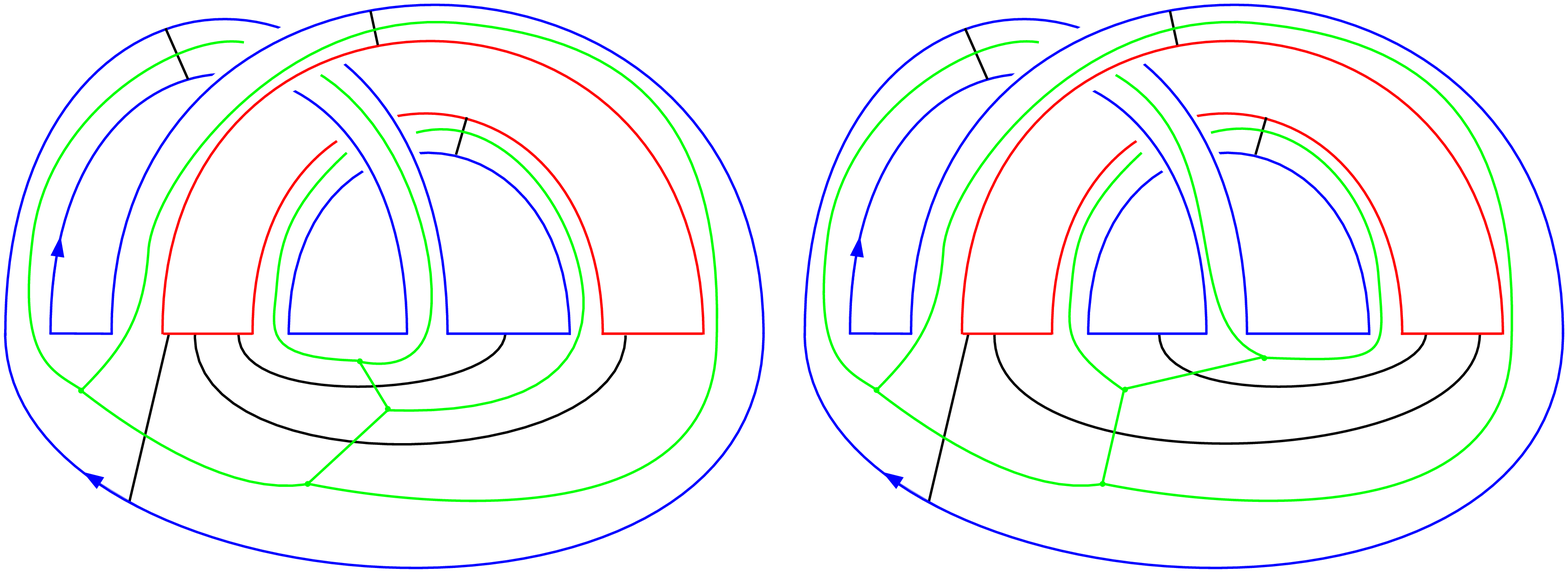}
	\caption{{\small  }}
	\label{fig: DWH-G2-3}
\end{figure}

When $\alpha$ splits $S_2$ into a connected genus one surface with two boundary components, we first need to generate all
possible $DWH$ diagrams along with all possible complete collections of weighted $\omega$-arcs.  Such a $DWH$ will
have three handles with each handle having the feature that it intersects one or both boundary curves.  Since we must
have at least one handle intersecting both boundary curves, the possibilities are:  all three handles intersect both boundaries
(as in Figure \ref{fig: DWH-G2-2}); two handles intersect both boundaries (as in Figure \ref{fig: DWH-G2-3})); or, only one handle intersects
both boundary curves.  In the latter case, it is straight forward to see that, due to the requirement that the number of $++$ and
$--$ arcs are equal, there will be boundary parallel $\omega$-arcs.  Since this would mean that $|\a \cap \b|$ is not
minimal, we conclude that only the first two possibilities occur.

Having settled on the $DWH$ diagram of either Figure \ref{fig: DWH-G2-2} or \ref{fig: DWH-G2-3}, we consider other choices for
a maximal collection of weighted $\omega$-arcs.  For example, Figure \ref{fig: DWH-G2-3}-left and \ref{fig: DWH-G2-3}-right
illustrates two different choices (w.r.t. the $DWH$ structure) in the case where we have just two handles
intersecting both boundaries.  Note that the difference
between the two collections is a different choice for the arcs associated with the weights $w_4, w_5, w_6$.  However,
by interchanging the roles of the $w_3$ and $w_4$---in Figure \ref{fig: DWH-G2-3}-right, we view $w_4$ as the co-core of a handle
instead of $w_3$---we obtain Figure \ref{fig: DWH-G2-3}-left (after a relabeling weights).  Moreover, there is a similar
re-choosing of the co-cores of handles in Figure \ref{fig: DWH-G2-3}-left that will yield the weighted arc collection of
Figure \ref{fig: DWH-G2-2}.  Finally, for the alternate choice of the $w_4, w_5, w_6$ $\omega$-arcs in Figure \ref{fig: DWH-G2-2},
one can again re-choose the handle co-cores to produce the collection of Figure \ref{fig: DWH-G2-2}.  (We leave the details to the reader.)
Thus, we need only consider the collection of weighted arcs in Figure \ref{fig: DWH-G2-2} and its associated ILP (\ref{LLP}).
\begin{figure}[htbp]
\labellist
\small\hair 1pt
\pinlabel $(9,10,2)$ at 226 140
\endlabellist
	\centering
		\includegraphics[width=0.75\textwidth]{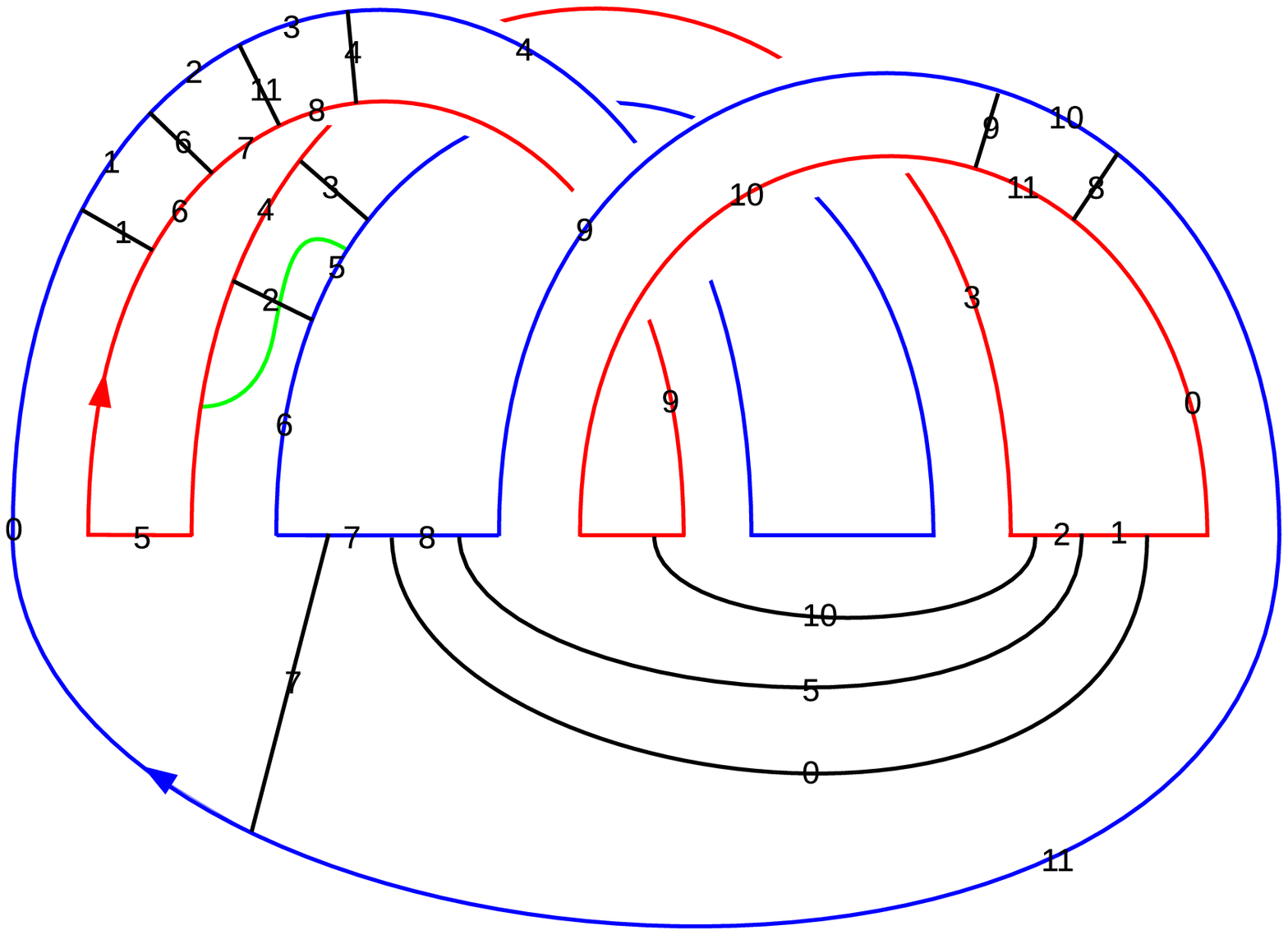}
	\caption{{\small }}
	\label{fig: DWH-G2-D4-I12}
\end{figure}
As mentioned at the end of \S\ref{sec: DWH}, the optimal value of ILP(\ref{LLP}) is $8$ and $[2,2,2,0,2,0]$ is the unique
weight solution.  Again, this was determined by utilizing the software package CPLEX \cite{[I]}.  Next, we employed MICC to determine that
any of the four associated filling pairs are only distance $3$.  Continuing, we utilized CPLEX to find all weight solutions for values $P=9,10,11$.
CPLEX found five weight solutions for value $9$: $[2, 3, 2, 0, 2, 0]$, $[3, 2, 2, 0, 2, 0]$, $[2, 2, 3, 0, 2, 0]$,
$[2, 2, 2, 0, 3, 0]$,and $[2, 2, 3, 0, 2, 0]$.
Using MICC we determined all associated intersection $P=9$ filling pairs are distance $3$.  Similarly, for values $P=10$ and $P=11$ there
are $55$ and $79$, respectively, weight solutions.  Again, MICC found only distance $3$ filling pairs.

However, using CPLEX to find all weight solutions to ILP (\ref{LLP}) for value $P=12$ we found $150$ solutions, $9$ of which had
associated filling pairs that are of distance $4$.  To list these solutions: $[w_1 , w_2 , w_3 , w_4 , w_5 , w_6] \in \{
[2, 2, 2, 2, 2, 2]^4,$ $[2, 4, 2, 0, 4, 0]^2,$ $[4, 2, 4, 0, 2, 0]^2,$ $[4, 0, 4, 2, 0, 2]^2,$ $[4, 2, 2, 1, 2, 1],$ $[2, 2, 4, 1, 2, 1],$ $[4, 0, 0, 2, 4, 2]^2,$ 
$[0, 4, 0, 2, 4, 2]^2,$ $[0, 4, 4, 2, 0, 2]^2\}$.  (The power notation is to dedicate a multiplicity of distance $4$ filling pairs.)
Figure \ref{fig: DWH-G2-D4-I12} is a realization of $[4, 2, 2, 1, 2, 1]$.  We credit its original discovery to J. Birman, D. Margalit and the
second author \cite{[BMM]}.  Finally, we remark that there are undoubtedly duplications up to homeomorphism which we do not
determine.

To summarize, utilizing the Distance $\geq 4$ Test capability of MICC, we have shown that the filling pair of
Figure \ref{fig: DWH-G2-D4-I12} is distance $\geq 4$.  To establish distance $4$, we must produce a length $4$ path
between $\alpha$ and $\beta$.  We refer the reader to the green arc in Figure \ref{fig: DWH-G2-D4-I12} that crosses the $\omega$-arc labeled $2$ and
has endpoints in the $5$ segments of $\partial_+$ and $\partial_-$.  Thus, this green arc can be understood as a closed curve that
intersects $\alpha$ once and $\beta$ once.  This closed curve, $\a_2$, will represent a vertex $v_2$ of a length $4$ path
$\{ v=v_0 , v_1 , v_2 , v_3 , v_4 = w\}$.  We can obtain a curve, $\a_1$, representing $v_1$ by taking a regular neighborhood of $\a_0 \cup \a_2$ and letting
$\a_1$ correspond to its unique boundary curve.  (Here, $\a = \a_0$ and $\a_4 = \b$.)  Similarly, $v_3$ is represented by a curve
coming from the boundary curve of a regular neighborhood
of $\a_4 \cup \a_2$. The fact that these neighborhoods are each topologically a torus-minus-disc makes all of these curves essential.
Thus, we have a length $4$ path.

\begin{figure}[htbp]
\labellist
\small\hair 2pt
\pinlabel $\g_1$ at 80 100
\pinlabel $\g_1$ at 470 90
\pinlabel $\g_1$ at 630 90
\pinlabel $\a_2$ at 265 90
\pinlabel $\g_2$ at 170 90
\pinlabel $\g_2$ at 522 90
\endlabellist
	\centering
		\includegraphics[width=0.80\textwidth]{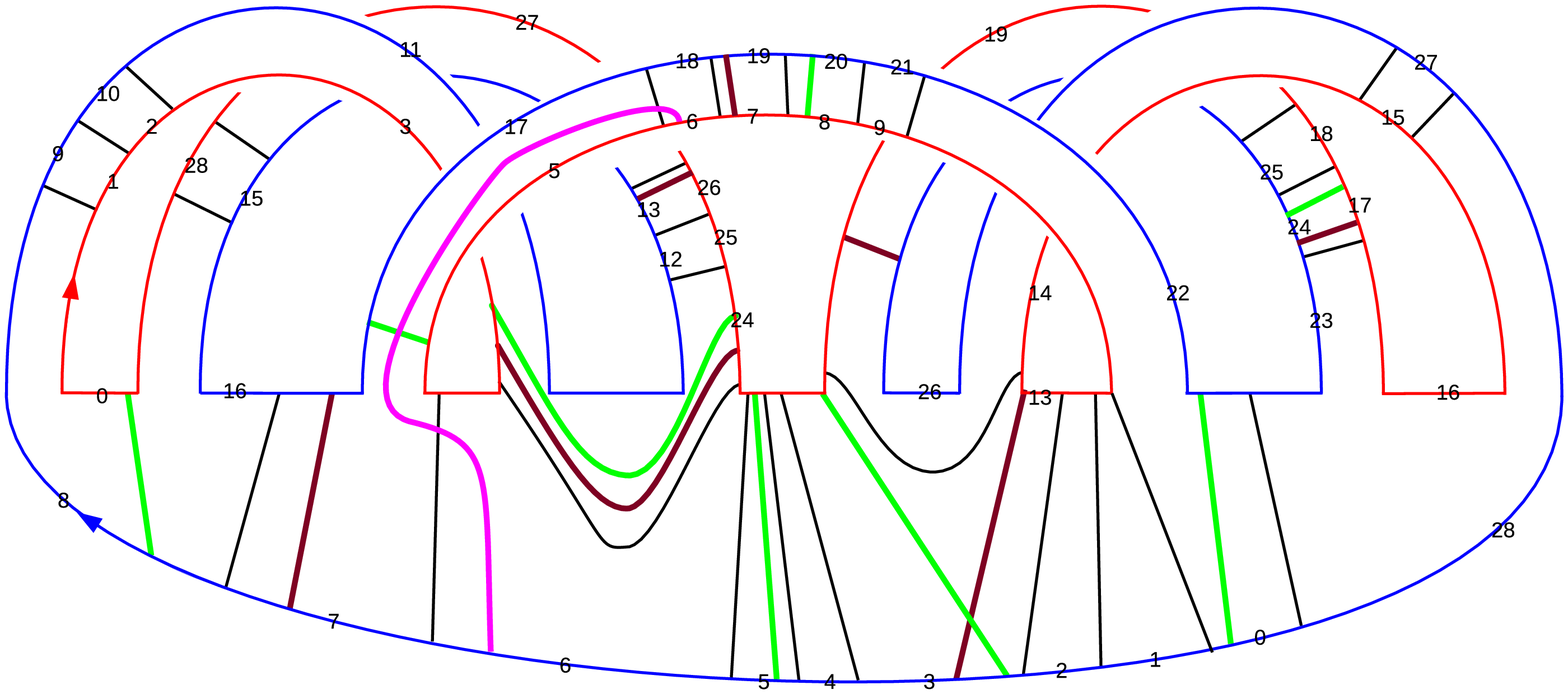}
	\caption{{\small }}
	\label{fig:genus 3}
\end{figure}

\subsection{Proof of Theorem \ref{theorem: DWH-G3-D4-I29}.}
We will prove the $DWH$ diagram of Figure \ref{fig:genus 3} is a genus $3$ distance $4$ example of a filling pair have intersection $29$.
As described in \S\ref{subsec: input}, the input for MICC can be obtained by traversing $\partial_+$ \& $\partial_-$ and reading
off the labels of the $\omega$-arcs as they are crossed.  The labels of the segments of the boundary curves indicate the
identification map that forms the curve $\alpha$.  As before, the black $\omega$-arcs join together to form $\b$.
For convenience, we give the MICC input below:
\begin{itemize}
\item[] {\tt Input top identifications: 1,11,3,27,8,15,7,24,0,10,2,12,4,21,19,17,\\ 
24,14,6,23,28,9,16,25,13,5,20,18,16}
\item[] {\tt Input bottom identifications: 0,10,2,26,7,14,6,23,28,9,1,11,3,22,20,\\
18,25,13,5,22,27,8,15,26,12,4,21,19,17}
\end{itemize}

In our genus $2$ example, we saw that there were $6$ curves representing vertices in the elementary circuit set $\G$.
In fact, MICC
utilizes the set $\G^\prime$ in its calculation.  For the filling pair of Figure \ref{fig:genus 3}, the set $\G^\prime$ has $28$ vertices.
We can specific a curve representing a vertex of $\G^\prime$ by a sequence of boundary segment labels.  For example, the
sequence {\tt [0, 22, 5, 17, 24, 3, 20, 8]} corresponds to a curve in this figure that intersects in cyclic order
the boundary segments in this list.  The set of green arcs in the figure correspond to this label sequence and should be to
understood as a curve, $\g_1$, representing $\bar\g_1 \in \G^\prime$.  By inspection observe that $(\a,\g_1)$ is a filling pair.
Similarly, the brown arcs of the figure correspond to the label sequence {\tt [3, 24, 17, 7, 19, 26, 13]}.  The reader can also
check by inspection that this curve, $\g_2$, representing $\bar\g_2 \in \G^\prime$, is also filling when paired with $\alpha$.

Having had MICC determine that all of the curves representing vertices in $\G^\prime$ fill when paired with $\alpha$, it remains for us to find a
length $4$ path so as to establish $d(v,w)=4$.  The magenta arc of Figure \ref{fig:genus 3} corresponds to
a curve, $\a_2$, representing a $v_2$ vertex
of such a path.  Notice $\a_2$ intersects $\b$ only twice.  Thus, as a pair they cannot be filling.  Taking an appropriate
boundary curve of a regular neighborhood of $\a_2 \cup \b$ will yield a $a_3$.  Similarly, since $|\a \cap \a_2|=1$ we
know $d(v,v_2)=2$ and we can construct an $\a_1$ as we did for our genus $2$ example to give us a representative
of $v_1$. \qed

\begin{remark}
\label{remark: experiment}
As previously mentioned, the argument establishing (the genus $2$) Theorem \ref{theorem: DWH-G2-D4-I12} can be thought of as a
proof of concept calculation.  For higher genus the calculation is exactly the same, although more extensive due to the fact that there are more $DWH$ representations and, thus, a $\Gamma$-calculation with accompanying ILP equations for each $DWH$.  However, if we willing to settle for an estimate the calculation can be limited to just one $DWH$ representation and we can use
MICC to ``discover'' some rough bounds.  This is essentially how the example of Figure \ref{fig:genus 3} was found.  We choose a somewhat
symmetric genus $3$ $DWH$ and placed an initial set of disjoint proper essential arcs so that any simple closed curve drawn in the $DWH$
intersected these arcs at least $5 (= [2 \cdot 3 - 1] =  [2 \cdot g - 1 ])$ times.
(For example, observe that a curve that goes over only the left most $1$-handle of Figure \ref{fig:genus 3} can be made to
intersect the black arcs exactly $5$ times.)
Using the {\tt perm} command discussed in \S\ref{subsec: perm}
we were able to find identifications of the two boundary curves that connected these arcs into a single curve.  From there we could use
MICC to compute the distance of the resulting filling pair.  When such an attempt yields a distance $3$ pair, we can use the
{\tt curves} command discussed in \S\ref{subsec: curves} to list all possible curves in $\Gamma$.  By adding in more proper arcs and/or
altering their placement and iterating the procedure above, we can reduce the size of $\Gamma$ until we find a distance $4$ filling pair.
Such experimentation may be useful in understanding the relationship between the distribution of the proper arcs in the $DWH$ and the
possible curves in $\Gamma$.
\end{remark}

\section{Information for the user of MICC.}
\label{sec: MICC commands}
In this section we discuss the input format for MICC and the commands for analyzing and manipulating that input.
The input is two sequences of numbers which corresponds to a ``ladder'' representation of curves in the surface.  The commands are {\em genus}, {\em distance}, {\em curves}, {\em matrix}, {\em faces} and {\em perm} (for permutation).  We will illustrate
each of these features by further developing the example in Figure \ref{fig: DWH-G2-D4-I12}.

\subsection{MICC input}
\label{subsec: input}
Given a $DWH$ presentation of filling pair of curves $(\a,\b)$ on a surface $S$, as previously discussed in \S\ref{sec: DWH},
we have specified the gluing of $\partial_+$ and $\partial_-$ by cyclically labeling their segments.  Assuming we used $k$ labels
for the segments of $\partial_\pm$, we now label the $\omega$-arcs of $\b$, $0$ through $k-1$.  Looking at Figure \ref{fig: DWH-G2-D4-I12}, we have a labeling using $0$ through $11$.  For the
purpose of extracting the MICC input from this example, the labeling assignment could have been done in a random fashion.
However, for aesthetic reasons we have taken care to label the $\omega$-arcs in the cyclic order they occur in $\b$.

With the $\omega$-arcs of $\b$ in our $DWH$ representation labeled as described above, we can now extract the MICC input for our filling pair.
Starting at segment $0$ in $\partial_-$, we traverse this boundary component in the positive direction and record
the labels of the $\omega$-arcs as we pass over their endpoints.  This sequence of labels will be the
input for MICC's ``{\tt Input top identifications}''.  Repeating this process on $\partial_+$ we get the sequence of labels that will
be the input for MICC's ``{\tt Input bottom identifications}''.  Together these two sequences are the basis of a {\em ladder representation},
$L_v(w)$, of a filling pair $(\a,\b)$.

\begin{figure}[htbp]
\labellist
\small\hair 2pt
\pinlabel $\a$ at 1 115
\pinlabel $\b$ at 25 155
\endlabellist

	\centering
		\includegraphics[width=0.80\textwidth]{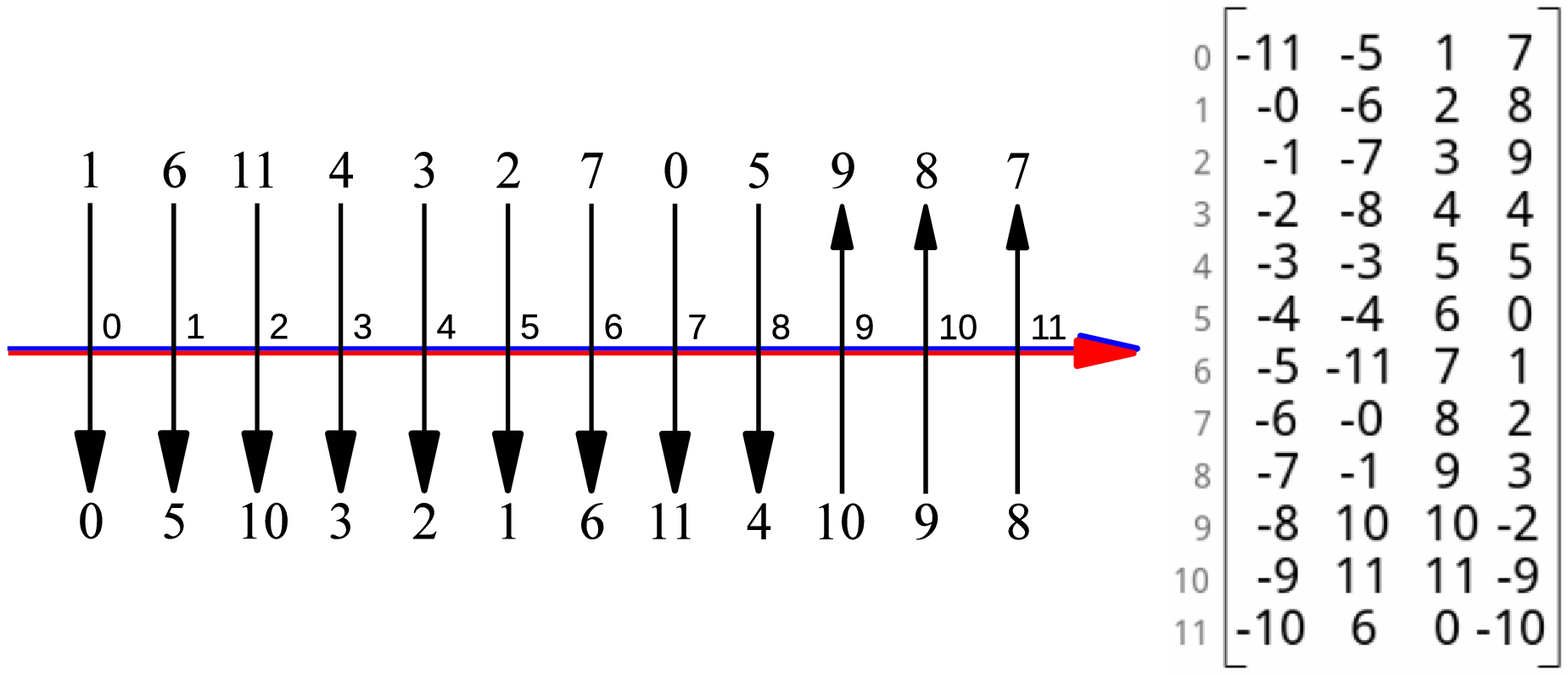}
	\caption{{\small }}
	\label{fig: d4-l-m}
\end{figure}

\begin{figure}[htbp]
\labellist
\small\hair 2pt
\pinlabel $\a$ at 7 75
\pinlabel $\b$ at 10 110
\pinlabel $v^-$ at 253 95
\pinlabel $w^+$ at 283 95
\pinlabel $v^+$ at 313 95
\pinlabel $w^-$ at 343 95
\endlabellist
	\centering
		\includegraphics[width=0.50\textwidth]{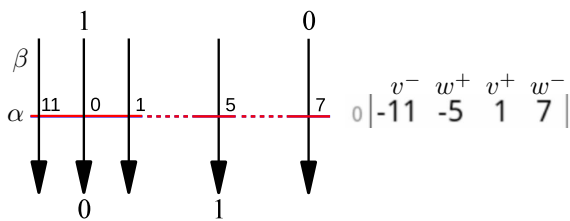}

	\caption{{\small The illustration depicts how to generate ${[L_v(w)]}_0 = [v^-(0), w^+(0) , v^+(0) , w^-(0)] = [-11, -5, 1 , 7]$.}}
	\label{fig: d4-l-to-m}
\end{figure}

A ladder representation is readily understood by again considering our example in Figure \ref{fig: DWH-G2-D4-I12}.
Starting at the $0$ segments of $\partial_\pm$ and reading off the two label sequences, our input for MICC would be:
\begin{itemize}
\item[] {\tt Input top identifications: 1,6,11,4,3,2,7,0,5,9,8,7}
\item[] {\tt Input bottom identifications: 0,5,10,3,2,1,6,11,4,10,9,8}
\end{itemize}
The reader should readily grasp the ladder metaphor by considering the representation
in Figure \ref{fig: d4-l-m}-left of a regular neighborhood of
$(\partial_+ \cup_\sim \partial_-) \cup \{{\rm all} \ \omega \ {\rm arcs} \}$ coming from Figure \ref{fig: DWH-G2-D4-I12}.
In this illustration, our two boundary components
which have been glued together to form $\alpha$ is represented by the horizontal segment that has its ends identified.
(Our convention forces this identification to always occur in the middle of the $0$-segment of $\partial_\pm$.)  Each 
vertical segment above or below this horizontal $\alpha$ is half of an $\omega$-arc.  From left-to-right there are $12$ labels above these
vertical $\omega$-halves correspond to MICC's {\tt Input top identifications}.  Similarly, the $12$ labels below these
vertical $\omega$-halves correspond to MICC's {\tt Input bottom identifications}.
It will be useful to the reader to have an understanding of how MICC utilizes a ladder representation.  Specifically, for a
ladder $L_v(w)$ we have, from left-to-right, $0$ through $k-1$ intersection points
of $\alpha$ with $\beta$.  For the $i^{\rm th}$ intersection we
will associate to it a $1 \times 4$ row vector ${[L_v(w)]}_i = [v^-(i),w^+(i),v^+(i), w^-(i)]$.  Figure \ref{fig: d4-l-to-m} illustrates the scheme
for determining the values of the $v^\pm {\rm 's}$ and $w^\pm {\rm 's}$.  In particular, $v^- (i)= -(i-1)$ and $v^+(i) = i+1$; and
$w^+(i)$ and $w^-(i)$ correspond to the adjacent (ladder) vertices along $\b$ with
the parity determined by whether $\b$ is pointing down ($w^+(i) <0$ and $w^-(i)>0$) or up ($w^+(i) >0$ and $w^-(i)<0$) at the $i^{\rm th}$ intersection.
Stacking the row vectors ${[L_v(w)]}_i$ in order, we produce a $ k \times 4$ matrix. Again, $k$ is the intersection number of $\alpha$ and $\beta$.
Figure \ref{fig: d4-l-m}-right is the associated {\em characteristic matrix}, $[L_v(w)]$,
for the ladder in Figure \ref{fig: d4-l-m}-left.  It should be readily apparent to the reader that the information in $[L_v(w)]$ is sufficient
to reproduce the ladder $L_v(w)$, and $L_v(w)$ is sufficient to reproduce a $DWH$ representation.  Thus, up to
any permutation of labels and changes of orientation of our curves, $L_v(w)$ is dependent only on the classes $v$ and $w$.

Given the prompt, {\tt What would you like to calculate?}, a reply of {\tt matrix} will produce $[L_v(w)]$.

\subsection{Genus command}
\label{sec: genus}
We now make concrete the filling calculation.
We will use $L_v(w)$ and its characteristic matrix $[L_v(w)]$ to compute the minimal genus of the surface that $\alpha$ and $\beta$ fill. For a surface of genus $g$, with a minimally intersecting filling $4$-valent graph, $\alpha \cup \beta$, we have
$g= -\frac{1}{2}(|V| - |E| + |F|)+1$, where $V = \a \cap \b$ is the set of vertices and $E$ is the set of edges of the graph $\a \cup \b$.  Also,
$F = S \setminus (\a \cup \b)$ is a set of $2n$-gon disc regions.
Then $|E| = 2|V|$.  Thus, determining $g$ requires a count of the number of $2n$-gon regions of $F$.

Listing the components of $F$ using $[L_v(w)]$ is achieved by specifying their edge-path boundaries in $\a \cup \b$.
When traversed such a edge-path will be {\em alternating}---an edge $e^{\a}_1 \in \a$ followed by an
edge $e^{\b}_1 \in \b$ followed by $ e^{\a}_2 \in \a$ and so on, will be a cyclically ordered set
$\{e^{\a}_1 , e^{\b}_1 , \cdots , e^{\a}_n, e^{\b}_n \}$.  Thus, $e^{\a}_j$ (respectively, $e^{\b}_j$) starts at vertex having integer label $\partial^s e^{\a}_j$
(respectively, $\partial^s e^{\b}_j$) and terminates at
vertex having integer label $\partial^t e^{\a}_j$ (respectively, $\partial^t e^{\b}_j$); and, $\partial^t e^{\b}_{j-1} = \partial^s e^{\a}_j$ and
$\partial^t e^{\a}_j = \partial^s e^{\b}_j$, all modulo $n$.

The entries of the ${[L_v(w)]_i}{\rm 's}$ gives us that
$$\partial^t e^{\a}_j \in \{ v^+(\partial^s e^{\a}_j) , |v^- (\partial^s e^{\a}_j)| \}  = \{ ( \partial^s e^{\a}_j -1 ) , (\partial^s e^{\a}_j +1 ) \}.$$
Additionally, we have
$$ \partial^s e^{\b}_j \in \{ |w^+ (\partial^t e^{\b}_j)| , |w^- (\partial^t e^{\b}_j)| \} \ {\rm and} \ 
\partial^t e^{\b}_j \in \{ |w^+ (\partial^s e^{\b}_j)| , |w^- (\partial^s e^{\b}_j)| \}.$$
Now fixing a $2n$-gon region and, by convention, traversing its
edge path boundary so as to always keep the region to our left, we have the following scheme for
finding the terminus endpoint for the {\em next} $e^\a$ or $e^\b$ edge.

\begin{itemize}
\item[T1--] If $\partial^t e^\a_j = v^+(\partial^s e^\a_j)$ then $\partial^t e^\b_j = |w^+ (\partial^t e^\a_j)|$ (with
$\partial^s e^\b_j = \partial^t e^\a_j$).
\item[T2--] If $\partial^t e^\a_j = |v^-(\partial^s e^\a_j)|$ then $\partial^t e^\b_j = |w^- (\partial^t e^\a_j)|$ (with
$\partial^s e^\b_j = \partial^t e^\a_j$).
\item[T3--] If $\partial^s e^\b_j = |w^+(\partial^t e^\b_j)|$ and $\partial^t e^\b_j = |w^- (\partial^s e^\b_j)|$
then $\partial^t e^\a_{j+1} = v^+(\partial^t e^\b_j)$ (with $\partial^s e^\a_{j+1} = \partial^t e^\b_j$).
\item[T4--] If $\partial^s e^\b_j = |w^-(\partial^t e^\b_j)|$ and $\partial^t e^\b_j = |w^+ (\partial^s e^\b_j)|$
then $\partial^t e^\a_{j+1} = |v^- (\partial^t e^\b_j)|$ (with $\partial^s e^\a_{j+1} = \partial^t e^\b_j$).
\item[T5--] If $\partial^s e^\b_j = |w^+(\partial^t e^\b_j)|$ and $\partial^t e^\b_j = |w^+ (\partial^s e^\b_j)|$
then $\partial^t e^\a_{j+1} = v^+ (\partial^t e^\b_j)$ (with $\partial^s e^\a_{j+1} = \partial^t e^\b_j$).
\item[T6--] If $\partial^s e^\b_j = |w^-(\partial^t e^\b_j)|$ and $\partial^t e^\b_j = |w^- (\partial^s e^\b_j)|$
then $\partial^t e^\a_{j+1} = v^- (\partial^t e^\b_j)$ (with $\partial^s e^\a_{j+1} = \partial^t e^\b_j$).
\end{itemize}

We illustrate this scheme using our $[L_v(w)]$ in Figure \ref{fig: d4-l-m}-right.  Starting at vertex $2$ we can
traverse the edge between $2(= \partial^s e^\a_1 = v^+ ( 1) )$ and $1 (= \partial^t e^\a_1 = |v^- (2)|)$.  Thus, we have the
assumption of T2 for $e^\a_1$.  This gives us
$\partial^t e^\b_1 = 8 (= |w^- (1 )| = \partial^s e^\a_2)$ with $1=\partial^s e^\b_1 = \partial^t e^\a_1$.  Since $1 = |w^+ (8)|$, we have the T3 assumption and
 $9 = \partial^t e^\a_2( = v^+(8))$.  But this gives us the T1 assumption for $e^\a_2$.
Thus, we have that
$\partial^t e^\b_2 = |w^+ (9)|$ which is $10$.  Since $|w^-(10)| = 9$ we again have the T3 assumption which implies $\partial^t e^\a_3 = |v^-(10)| = 9$.
Now, we again have the T2 assumption for $e^\a_3$ which implies $\partial^t e^\b_3 = | w^- (9)| = 2$, back where we started.  (Refer to the region in
Figure~\ref{fig: DWH-G2-D4-I12} contain $(9,10,2)$.)

Working out all such boundary edge paths we get a count for $|F|$ and, thus, are able to compute genus.

\subsection{Faces command}
\label{subsec: faces}
Having traversed all the boundaries of the $2n$-gon regions of $S \setminus (\a \cup \b)$, MICC records this
calculation as a vector.  Specifically, let $ F_{2n} \subset F$ be the number of $2n$-gon regions for $n \in \{2,3,4, \cdots \}$.
Then associated with the graph $\a \cup \b \subset S_g$ we have the vector $[F_4 , F_6 , F_8 , \cdots]$.
The output of MICC is actually formatted as $\{4:\ F_4 , \  6:\ F_6, 8: \ F_8, \cdots \}$.  This
vector solution to Euler characteristic equation of Lemma 4.1 of \cite{[He1]}.
The prompting inquiry and command appears as:  {\tt What would you like to calculate? faces}.  
Additionally, it lists each boundary edge in a truncated matter---it lists only the $e^\a_i {\rm 's}$.
For our extended example of
Figure \ref{fig: d4-l-to-m}, we would get:
\begin{eqnarray}
\begin{tabular}{l}
${\tt Vector solution:}  \  \{4: 6, 6: 4\}$ \\
$(0, 11, 7)$ \\
$(0, 5, 6)$ \\
$(1, 6)$ \\
$(8, 1)$ \\
$(2, 7)$ \\
$(9, 10, 2)$ \\
$(8, 3)$ \\
$(9, 3, 4)$ \\
$(4, 5)$ \\
$(10, 11)$ \\
\end{tabular}
\end{eqnarray}
\subsection{Curves command}
\label{subsec: curves}
An alternate reply to the prompt, {\tt What would you like to calculate?}, is {\tt curves}.  MICC applies the Theorem \ref{theorem: BMM test} test by listing
all of the curves representing vertices of $\G^\prime$ and computing the genus of their graphs when paired with $\a$.
MICC finds all curves representing elements of $\G^\prime$ by applying a classical depth-first search \cite{[Si]} for elementary circuits of the
graph $G(C^\prime)$, the dual graph to the proper arcs of $\a$ in $S \setminus \b$.

The MICC output is a cyclic sequence of $e^{\a}$ edges.  That is, a curve,$\g$, representing a vertex
$\bar\g \in \G^\prime$ and a region $f \in F$, $\g \cap f$ will be a
collection of proper arcs having their endpoints on $e^{\a}$ edges of the boundary of $f$.  As we traverse $\g$, we will travel between
regions by passing through $e^{\a}{\rm's}$.  Thus, $\bar\g$ can be characterized by giving a cyclic listing of these $e^{\a}{\rm's}$.

Continuing with our extended example of Figure \ref{fig: d4-l-to-m}, the output response to {\tt curves} would yield:
\begin{eqnarray}
\begin{tabular}{l}
${\tt Path \  [0, 7, 2, 9, 3, 8, 1, 6]}$ \\
${\tt Curve genus:} \  2$ \\ \\
${\tt Path \  [2, 10, 11, 7]}$ \\
${\tt Curve genus:} \  2$ \\ \\
${\tt Path \  [1, 6, 5, 4, 3, 8]}$ \\
${\tt Curve genus:} \  2$ \\ \\
${\tt Path \  [0, 5, 4, 9, 2, 7]}$ \\
${\tt Curve genus:} \  2$ \\ \\
${\tt Path \  [0, 5, 4, 9, 10, 11]}$ \\
${\tt Curve genus:} \  2$ \\ \\
${\tt Path \  [0, 11, 10, 9, 3, 8, 1, 6]}$ \\
${\tt Curve genus:} \  2$ \\
\end{tabular}
\end{eqnarray}

\subsection{Distance command}
\label{subsec: distance}
To determine $d(v,w)$ one replies to the prompt, {\tt What would you like to calculate?}, by typing {\tt distance}.
If the genus of any of the pairs $(\a,\g)$ is less than $g$, then MICC will respond with {\tt Distance:  3}.
If the genus of all such pairs is $g$, then MICC will respond with {\tt Distance:  4+}.
All possible $\g$ representing elements of $\G^\prime$ are determined as described above.

\subsection{Perm command}
\label{subsec: perm}
Finally, MICC has an experimental functionality.  As most topologist who have attempted to construct filling pairs on
surfaces know, it is difficult to do so while avoiding the production of multi-curves.  For example, starting with a $DWH$ representation
of the surface, after placing down some collection of $\omega$-arcs finding an identification of $\partial^+$ and $\partial^-$
so as to have $\b$ be a single curve is tedious at best.  Fortunately, MICC automates this process.  Given any
ladder top/bottom identification, it will first determine whether $\b$ is a single curve or a multi-curve.  If it is a multi-curve, then
it will produce the inquiring prompt  {\tt Would you like to shear this multi-curve?}.  With the reply {\tt yes}, it will search through
all possible $\partial^+$/$\partial^-$ identifications for those that yield a single $\b$ curve and print them out along with
their distance.  If MICC has been given a ladder identification that corresponds to a single $\b$ curve, one can still
find all other $\partial^+$/$\partial^-$ identifications that yield a single curve.  When given the prompt 
{\tt What would you like to calculate?} just reply with {\tt perm} (for permutation).

For our extended example, the output would be:
\begin{eqnarray}
\begin{tabular}{l}
${\tt Curve \  1 \  Distance: \   3}$ \\
$[2, 7, 12, 5, 9, 8, 7, 1, 6, 11, 4, 3]$ \\
$[12, 5, 10, 3, 2, 1, 6, 11, 4, 10, 9, 8]$ \\ \\
${\tt Curve \ 2 \ Distance: \ 3}$ \\
$[5, 9, 8, 7, 1, 6, 11, 4, 3, 2, 7, 12]$ \\
$[12, 5, 10, 3, 2, 1, 6, 11, 4, 10, 9, 8]$ \\ \\
${\tt Curve \ 3 \ Distance:\  4+}$ \\
$[1, 6, 11, 4, 3, 2, 7, 12, 5, 9, 8, 7]$ \\
$[12, 5, 10, 3, 2, 1, 6, 11, 4, 10, 9, 8]$
\end{tabular}
\end{eqnarray}

\section{Concluding remarks.}
\label{sec: conclusion}
MICC was originally created as a tool that would help in the search for distance $4$
filling pairs on surface of genus greater than $2$. We have since realized that it can be set to other uses.
For example, it was recently used to find geodesic triangles where any pair of vertices correspond to filling pairs in minimal position.  As mentioned previously, we hope researchers will find additional uses.

Remarking on the complexity of the algorithms employed in MICC, the most computationally expensive task involves finding all cycles in the graph $G(A^\prime)$
of \S \ref{subsec: curves}. We know that the runtime of the process of counting all cycles in a graph is at best exponential.
This gives a lower bound on the complexity of finding all cycles in $G(A^\prime)$, and thus on the complexity of the program itself.  

There are many ways in which we would like to improve the current version of MICC. Currently, MICC is only a partial implementation of the algorithm presented in \cite{[BMM]}. We plan to extend the functionality of MICC to encompass the full scope of the Efficient Geodesic Algorithm of \cite{[BMM]}. Yet this partial implementation is also manifestation of the complexity barrier involved with the current exponential running time of the graph search. An improved algorithm will allow for more intricate curves to be studied, and parallelization of the MICC would help future users to fully utilize their multicore computers in their research.

\noindent
{\small{\bf Acknowledgements}} The first, third and fourth authors are grateful to Joan Birman, Dan Margalit and the second author for
sharing results of their joint work as it developed.
Our thanks goes to John Ringland, Joaquin Carbonara and the URGE to Compute
program at the University at Buffalo and Buffalo State College for supplying a nurturing environment for our work.
This work was supported in part by NSF CSUMS grants 
0802994 and 0802964 to the University at Buffalo and Buffalo State College.
Finally, we thank the referees for alerting us to the recent results
in the literature, for suggesting a strengthening of the statement of Corollary \ref{corollary: DWH-G2-D4-I12} and numerous other expository improvements.

\section{Addendum:  All weight solutions of distance $4$ filling pairs in $\mathcal{C}^1(S_2)$ for $12$, $13$ \& $14$ intersections
when one curve is non-separating.}
Below we list with multiplicity all solutions to ILP (\ref{LLP}) for $P$-values $12$, $13$ and $14$.

For $P=w_1 + w_2 + w_3 + w_4 + w_5 + w_6 = 12$:
\begin{itemize}
\item[] $[2, 2, 2, 2, 2, 2]^4$, $[2, 4, 2, 0, 4, 0]^2$, $[4, 2, 4, 0, 2, 0]^2$, $[4, 0, 4, 2, 0, 2]^2$,
$[4, 2, 2, 1, 2, 1]$, $[2, 2, 4, 1, 2, 1]$, $[4, 0, 0, 2, 4, 2]^2$, $[0, 4, 0, 2, 4, 2]^2$, $[0, 4, 4, 2, 0, 2]^2$
\end{itemize}

For $P=w_1 + w_2 + w_3 + w_4 + w_5 + w_6 = 13$:
\begin{itemize}
\item[] $[5, 4, 2, 0, 2, 0]$, $[2, 3, 2, 2, 2, 2]^2$, $[2, 2, 5, 1, 2, 1]$, $[4, 1, 2, 1, 4, 1]$, $[4, 5, 2, 0, 2, 0]$,
$[1, 4, 4, 1, 2, 1]$, $[1, 4, 4, 2, 0, 2]$, $[2, 2, 3, 1, 4, 1]$, $[5, 2, 2, 1, 2, 1]$, $[5, 2, 2, 0, 4, 0]$,
$[2, 4, 5, 0, 2, 0]$, $[0, 4, 1, 2, 4, 2]$, $[2, 4, 4, 1, 1, 1]$, $[3, 2, 2, 2, 2, 2]$, $[4, 2, 2, 0, 5, 0]$,
$[2, 4, 3, 0, 4, 0]^2$, $[3, 4, 2, 0, 4, 0]^2$, $[4, 2, 4, 0, 3, 0]^2$, $[4, 0, 4, 2, 1, 2]$, $[4, 3, 2, 1, 2, 1]$,
$[4, 1, 0, 2, 4, 2]$, $[2, 4, 3, 1, 2, 1]^2$, $[2, 2, 2, 1, 5, 1]$, $[2, 2, 4, 1, 3, 1]$, $[4, 1, 2, 2, 2, 2]$,
$[2, 3, 4, 1, 2, 1]^2$, $[3, 4, 2, 1, 2, 1]$, $[2, 2, 5, 0, 4, 0]$, $[2, 5, 4, 0, 2, 0]$, $[2, 2, 2, 2, 3, 2]$
$[4, 0, 1, 2, 4, 2]$, $[3, 2, 2, 1, 4, 1]^2$, $[2, 2, 4, 0, 5, 0]$, $[1, 4, 0, 2, 4, 2]$, $[2, 5, 2, 1, 2, 1]$,
$[1, 4, 2, 2, 2, 2]$, $[2, 2, 1, 2, 4, 2]$, $[2, 2, 3, 2, 2, 2]^2$, $[0, 4, 4, 2, 1, 2]$, $[4, 2, 1, 1, 4, 1]$,
$[4, 1, 4, 2, 0, 2]$, $[2, 2, 4, 2, 1, 2]$, $[4, 3, 4, 0, 2, 0]^2$, $[4, 2, 2, 1, 3, 1]^2$
\end{itemize}

For $P=w_1 + w_2 + w_3 + w_4 + w_5 + w_6 = 14$:
\begin{itemize}
\item[] $[4, 0, 6, 1, 2, 1]^2$, $[3, 2, 4, 2, 1, 2]$, $[4, 0, 2, 3, 2, 3]^2$, $[0, 4, 4, 2, 2, 2]$, $[0, 4, 2, 1, 6, 1]^2$,
$[4, 3, 1, 1, 4, 1]$, $[4, 2, 2, 1, 4, 1]^4$, $[0, 4, 2, 2, 4, 2]$, $[2, 2, 6, 0, 4, 0]$, $[1, 4, 2, 2, 3, 2]$
$[4, 2, 4, 2, 0, 2]$, $[5, 3, 2, 0, 4, 0]$, $[2, 2, 4, 2, 2, 2]^3$, $[3, 3, 2, 1, 4, 1]$,
$[6, 2, 0, 1, 4, 1]^2$, $[4, 3, 4, 0, 3, 0]^2$, $[2, 4, 2, 2, 2, 2]^3$, $[2, 6, 4, 0, 2, 0]$, $[4, 4, 4, 0, 2, 0]^2$,
$[2, 3, 2, 2, 3, 2]$, $[2, 2, 5, 1, 3, 1]$, $[6, 2, 4, 1, 0, 1]^2$, $[2, 2, 6, 1, 2, 1]$, $[4, 4, 2, 1, 2, 1]^2$,
$[2, 2, 2, 3, 2, 3]^2$, $[3, 2, 3, 2, 2, 2]$, $[3, 5, 2, 1, 2, 1]$, $[4, 1, 4, 2, 1, 2]$, $[3, 5, 4, 0, 2, 0]$,
$[6, 2, 2, 0, 4, 0]$, $[4, 2, 3, 0, 5, 0]$, $[3, 3, 2, 2, 2, 2]^2$, $[4, 1, 3, 2, 2, 2]$, $[2, 2, 3, 2, 3, 2]^2$,
$[1, 4, 4, 1, 3, 1]$, $[6, 2, 4, 0, 2, 0]^2$, $[3, 3, 4, 1, 2, 1]$, $[2, 3, 5, 0, 4, 0]$, $[4, 0, 2, 2, 4, 2]$,
$[4, 2, 4, 1, 2, 1]$, $[3, 4, 4, 1, 1, 1]$, $[2, 4, 2, 0, 6, 0]^2$, $[0, 4, 6, 1, 2, 1]^2$, $[2, 2, 3, 1, 5, 1]$,
$[1, 4, 1, 2, 4, 2]$, $[4, 6, 2, 0, 2, 0]$, $[2, 4, 6, 0, 2, 0]$, $[2, 4, 4, 0, 4, 0]^2$, $[2, 4, 4, 1, 2, 1]^4$,
$[4, 2, 2, 0, 6, 0]$, $[2, 2, 0, 3, 4, 3]^2$, $[2, 2, 4, 3, 0, 3]^2$, $[4, 1, 3, 1, 4, 1]$, $[2, 4, 4, 2, 0, 2]$,
$[3, 4, 3, 1, 2, 1]$, $[2, 6, 2, 1, 2, 1]$, $[3, 2, 3, 1, 4, 1]$, $[4, 2, 0, 2, 4, 2]$, $[4, 0, 4, 2, 2, 2]$,
$[2, 2, 4, 1, 4, 1]^2$, $[4, 2, 2, 2, 2, 2]^3$, $[2, 3, 1, 2, 4, 2]$, $[0, 4, 2, 3, 2, 3]^2$, $[2, 4, 5, 0, 3, 0]$
$[5, 3, 2, 1, 2, 1]$, $[2, 4, 3, 1, 3, 1]$, $[4, 2, 4, 0, 4, 0]^2$, $[4, 3, 2, 1, 3, 1]$, $[4, 4, 2, 0, 4, 0]^2$,
$[2, 2, 2, 1, 6, 1]$, $[2, 3, 4, 1, 3, 1]$, $[6, 4, 2, 0, 2, 0]$, $[5, 4, 2, 0, 3, 0]$, $[4, 0, 2, 1, 6, 1]^2$,
$[4, 2, 6, 0, 2, 0]^2$, $[4, 5, 3, 0, 2, 0]$ $[2, 4, 2, 1, 4, 1]$, $[4, 2, 3, 1, 3, 1]$, $[6, 2, 2, 1, 2, 1]$,
$[2, 2, 2, 2, 4, 2]^3$, $[2, 6, 2, 0, 4, 0]^2$, $[3, 2, 4, 0, 5, 0]$, $[3, 4, 3, 0, 4, 0]^2$, $[2, 6, 4, 1, 0, 1]^2$,
$[2, 2, 4, 0, 6, 0]$, $[2, 6, 0, 1, 4, 1]^2$.
\end{itemize}

\begin{figure}[htbp]
	\centering
		\includegraphics[width=0.80\textwidth]{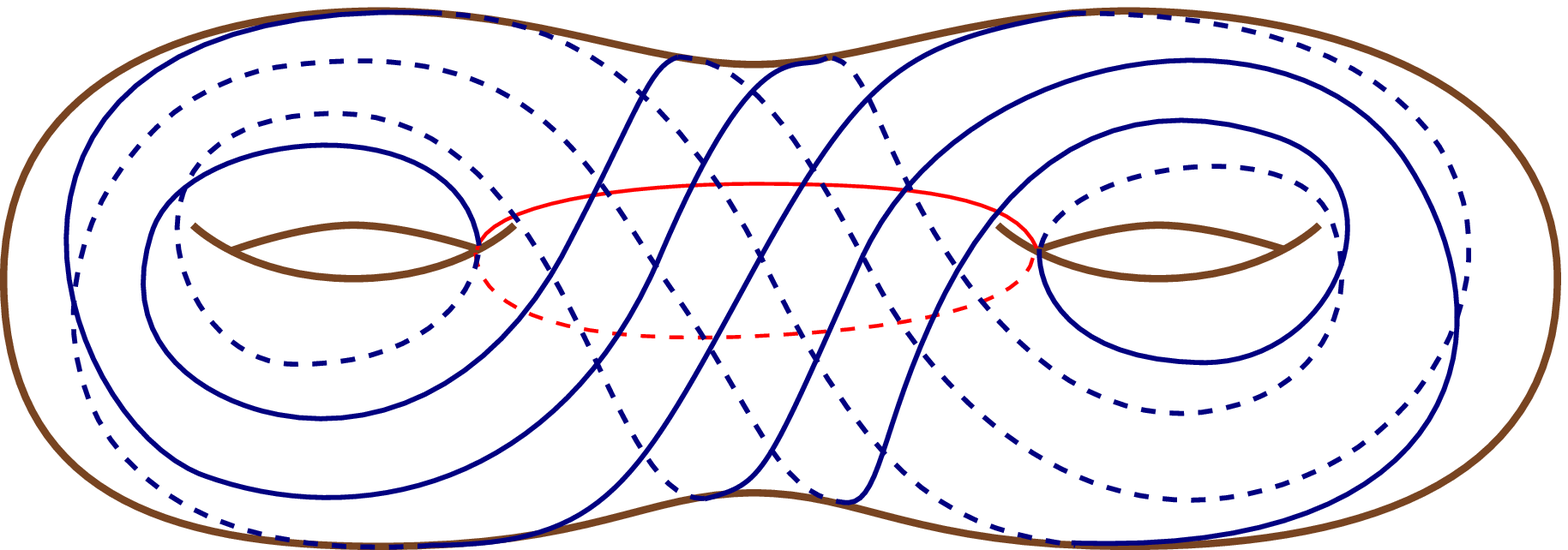}
	\caption{{\small }}
	\label{fig: 2,2,2,2,2,2}
\end{figure}

The weight solutions $[2,2,2,2,2,2]$ for $P=12$ is particularly intriguing since it suggests there might be a high level of symmetry.
Figure \ref{fig: 2,2,2,2,2,2} is a $3$-dimensional rendering of the one the four associated distance $4$ filling pair.  The aesthetic of this
rendering is so appealing that it was placed at the begin of \cite{[BMM]}.  We thank those authors for its use here.

MICC software package, software tutorial, and all known weight solutions yielding filling pairs having distance $\geq 4$ for
$P \leq 25$ with $g=2$ (approximately 72,000 weight solutions) are posted for download at micc.github.io.

\end{document}